\documentclass[twoside,magyar,english]{article}
\usepackage[T1]{fontenc}
\usepackage[latin2,latin1]{inputenc}
\usepackage{amsmath}
\usepackage{amssymb}

\makeatletter
\date{}

\AtBeginDocument{

}

\usepackage{babel}
\makeatother
\begin{document}

\newcommand{\zn}[1]{\left(\mathbb{Z}/#1\mathbb{Z}\right)^{*}}

\newcommand{\SL}{SL_{2}\left(\mathbb{Z}\right)}

\newcommand{\id}{\mathbb{I}}

\newcommand{\pk}{\mathbb{K}}

\newcommand{\set}[2]{\left\{  #1\, |\, #2\right\}  }

\newcommand{\cyc}[1]{\mathbb{Q}\left[ \zeta_{#1 }\right] }

\newcommand{\Mod}[3]{#1\equiv#2\, \left(\mathrm{mod}\, \, #3\right)}

\newcommand{\FA}[1]{\Gamma\left(K,\mathfrak{h}\right)}

\newcommand{\FB}[1]{\mathfrak{V}}

\newcommand{\FC}[1]{\Gamma}

\newcommand{\FD}[1]{\mathcal{S}_{g}}

\newcommand{\FE}[2]{\mathcal{M}_{\phi}}

\newcommand{\FF}[2]{}

\title{Mapping class group representations and Conformal Field Theory}

\author{P. Bantay}

\maketitle
\begin{center}Institute for Theoretical Physics, Eötvös University,
Budapest\end{center}

\begin{abstract}
We discuss some properties of the tower of mapping class group representations
associated to a Rational Conformal Field Theory. In particular, after
reviewing the elementary properties of the modular representation,
we discuss the Galois action, the structure of the projective kernel,
and the trace identities generalizing the formula of Verlinde.
\end{abstract}

\section{Introduction}

String theory and the closely related Conformal Field Theories have
led to many exciting interactions between mathematics and theoretical
physics in the past two decades. Among other things, they stimulated
the creation of two new mathematical concepts: that of Vertex (Operator)
Algebras \cite{Bo,FLM} and of Modular Tensor Categories \cite{Turaev,MTC}.
Both have found important applications in the sequel, e.g. Vertex
Algebras are the basis of our present understanding of Moonshine (the
relation between the representation theory of the Monster - the largest
sporadic simple group - and modular functions), while Modular Tensor
Categories provide new topological invariants of 3-manifolds \cite{Witten,Turaev}.
While Vertex Operator Algebras may be understood as the mathematical
formulation of the chiral algebras of Conformal Field Theory, Modular
Tensor Categories formalize the so-called Moore-Seiberg data \cite{MS1,MS}.
The above two concepts are closely related, as follows from Huang's
recent proof of the Verlinde conjecture \cite{Huang}.

A Modular Tensor Category provides, among other things, a sequence
\[
\FB{}_{1},\FB{}_{2},\ldots\]
of finite dimensional linear spaces that afford representations of
the mapping class groups\[
\FC{}_{1},\FC{}_{2},\ldots\]
of closed surfaces, where $\FC{}_{g}$ denotes the mapping class group
of a surface of genus $g$. This sequence of representations, called
sometimes the ''modular tower'', is very special, with many relations
between its different terms, and it does characterize the Modular
Tensor Category to a great extent. The aim of the present paper is
to review some properties of this sequence that could prove useful
in attempts to classify Modular Tensor Categories. As one may expect,
the best understood term is the first one, which is nothing but a
finite dimensional representation of the classical modular group $\FC{}_{1}=\SL$,
but we'll see that this $\SL$ module is already very special in many
respects, and does have a profound influence on the structure of the
subsequent terms.

In the next section we'll summarize some properties of the $\SL$
module $\FB{}_{1}$ which form the basis of most results to follow.
Section \ref{sec:The-Galois-action} describes what we know about
the number theoretic properties of this representation, while Section
\ref{sec:The-projective-kernel} discusses the structure of its (projective)
kernel. In Section \ref{sec:Trace-identities} we'll turn to the study
of the higher terms of the modular tower, and see the deep influence
the first term $\FB{}_{1}$ has on their structure. Finally, we'll
discuss briefly the implications of these results on the problem of
enumerating systematically Modular Tensor Categories.

\section{The modular representation\label{sec:The-modular-representation}}

As alluded to in the introduction, each Modular Tensor Category determines
a finite dimensional module $\FB{}_{1}$ of the classical modular
group $\FC{}_{1}=\SL$. Actually much more is true, for this module
comes equipped with a distinguished basis whose elements correspond
to the simple objects of the category, and the representation matrices
relative to this basis enjoy some remarkable properties. To explain
these properties, let's first recall that the group $\SL$, consisting
of the unimodular two-by-two integer matrices, is generated by the
matrices $\left(\begin{array}{cc}
1 & 1\\
0 & 1\end{array}\right)$ and $\left(\begin{array}{cc}
0 & -1\\
1 & 0\end{array}\right)$. As is customary, we'll denote by $T$ and $S$ respectively the
matrices representing these two $\SL$ elements in the above mentioned
distinguished basis. Of course, the matrices $T$ and $S$ have to
satisfy the defining relations of the group $\SL$, namely $S^{4}=1$
and the modular relation \begin{equation}
STS=T^{-1}ST^{-1}\,.\label{eq:modrel}\end{equation}

The fundamental properties of the modular representation are \cite{Verlinde,MS}:

\begin{enumerate}
\item The matrix $T$ is diagonal and has finite order;
\item The matrix $S$ is symmetric;
\item $S^{2}$ is a permutation matrix (of order 2);
\item There is a row of $S$ (usually labeled by $0$) all of whose entries
are positive;
\item The expressions\begin{equation}
\sum_{s}\frac{S_{ps}S_{qs}S_{rs}}{S_{0s}}\,,\label{eq:Verlinde}\end{equation}
related to the decomposition of tensor products of simple objects,
yield non-negative integers for any choice of simple objects $p,q,r$
(remember that the rows of $S$ are labeled by the simple objects
of the category).
\end{enumerate}
All these properties are satisfied by the modular representation associated
to a Modular Tensor Category, but it is not too hard to construct
examples of pair of matrices $T$ and $S$ that satisfy all of them,
but do not correspond to a Modular Tensor Category, i.e. requiring
them is a necessary but not sufficient condition. This raises the
question whether it is possible to give a set of sufficient conditions,
a problem that hasn't been settled yet. Extra necessary conditions
arise from the theory of the Galois action and the trace identities
to be explained later, but there are some others that we do know about,
e.g. those related to the existence of Frobenius-Schur indicators.
These later arise because to each simple object $p$ of a Modular
Tensor Category one can associate an analog of the classical Frobenius-Schur
indicator, i.e. a quantity $\nu_{p}$ that can take on the three values
$-1,0$ and $+1$, and which characterizes the symmetry of the tensor
square of $p$. The crucial thing is that this indicator is fully
determined by the modular representation, for one has the formula
\cite{FS}\begin{equation}
\nu_{p}=\sum_{q}\frac{S_{pq}\left|M(2)_{0q}\right|^{2}}{S_{0q}}\,,\label{eq:FS}\end{equation}
where $M(k)=S^{-1}T^{k}S$ for integer $k$. That the rhs. of Eq.(\ref{eq:FS})
can only be $-1,0$ or $+1$ does not follow from the properties of
the matrices $T$ and $S$ explained above, thus it provides a set
of further necessary conditions for the modular representation to
correspond to some Modular Tensor Category.

This is far from being the end of the story. Having a representation
of the group $\SL$, it is natural to ask whether its image is finite
or not, i.e. if $T$ and $S$ generate a finite matrix group. As it
turns out, this is always true, and even much more. To explain the
precise result, let's recall that for a positive integer $n$, the
principal congruence subgroup $\Gamma\left(n\right)$ of level $n$
consists of those elements of $\SL$ that are congruent to the identity
matrix modulo $n$, i.e. 

\begin{equation}
\Gamma\left(n\right)=\set{\left(\begin{array}{cc}
a & b\\
c & d\end{array}\right)\in\SL}{\Mod{a,d}{1}{n}\,\,\mathrm{and}\,\,\Mod{b,c}{0}{n}}\,.\label{eq:pcongdef}\end{equation}
Clearly, one has $\Gamma\left(1\right)=\SL$, and $\Gamma\left(n\right)$
is always a normal subgroup of $\SL$, being the kernel of the natural
homomorphism $\SL\rightarrow SL_{2}\left(\mathbb{Z}/n\mathbb{Z}\right)$.
A subgroup of $\SL$ is called a congruence subgroup of level $n$
if it contains $\Gamma\left(n\right)$, but no $\Gamma\left(k\right)$
for $k<n$.

The basic result about the kernel is that it is a congruence subgroup
of level $N$, where $N$, the so-called conductor, is equal to the
order of the matrix $T$ (remember that this order is always finite)
\cite{CG2,CMP}. As $\Gamma\left(N\right)$ has finite index in $\SL$,
this implies that $T$ and $S$ generate a finite matrix group, a
homomorphic image of $SL_{2}\left(\mathbb{Z}/NZ\right)$.

An obvious question in this context is whether the conductor $N$
could be arbitrary. One may show, using the theory of the Galois action
to be explained in the next section, that $N$ is bounded from above
by a function of the number of simple objects, i.e. the dimension
$\dim\FB{}_{1}$ of the modular representation \cite{CMP}. This upper
bound plays a crucial role in attempts to enumerate systematically
Modular Tensor Categories.

\section{The Galois action\label{sec:The-Galois-action}}

As should be clear from the previous section, the modular representation
associated to a Modular Tensor Category has indeed very special features.
In particular, because the simple objects provide a distinguished
basis in $\FB{}_{1}$, the matrix elements of $T$ and $S$ have an
invariant meaning, and one might ask about their arithmetic properties.
This is the subject of the theory of the Galois action \cite{BG,CG1}.

The basic idea is to look at the field $F$ obtained by adjoining
to the rationals $\mathbb{Q}$ all matrix elements of $T$ and $S$.
Note that, because $T$ and $S$ generate the image, the representation
matrix of any element of $\SL$ will have matrix elements lying in
$F$. From the properties of the modular representation explained
in the previous section it follows that $F$ is a finite Galois extension
of $\mathbb{Q}$, whose Galois group is abelian. By the famous theorem
of Kronecker-Weber, this means that $F$ is a subfield of some cyclotomic
extension of $\mathbb{Q}$, and one may even show that $F=\cyc{N}$,
where $N$ denotes the conductor and $\zeta_{N}$ is a primitive $N$-th
root of unity. By known results of algebraic number theory, this means
that the Galois group $Gal\left(F/\mathbb{Q}\right)$ is isomorphic
to the group $\zn{N}$ of prime residues mod $N$, the action of $\sigma_{l}\in Gal\left(F/\mathbb{Q}\right)$
corresponding to $l\in\zn{N}$ being determined by \begin{equation}
\sigma_{l}\left(\zeta_{N}\right)=\zeta_{N}^{l}\,.\label{eq:sigmadef}\end{equation}

As all matrix elements of $T$ and $S$ belong to $F$ by definition,
one may ask which matrix one gets by applying $\sigma_{l}$ to $T$
and $S$ element-wise. In case of $T$ the result is simply \begin{equation}
\sigma_{l}\left(T\right)=T^{l}\,,\label{eq:Ttrans}\end{equation}
because $T$ is diagonal, and its diagonal entries are all powers
of $\zeta_{N}$. In case of $S$ the result reads \cite{CG1}\begin{equation}
\sigma_{l}\left(S\right)=SG_{l}\,,\label{eq:Strans}\end{equation}
where the orthogonal and monomial matrices $G_{l}$ form a representation
of the group $\zn{N}$, i.e. \begin{equation}
G_{lm}=G_{l}G_{m}\,.\label{eq:Gmult}\end{equation}
Monomiality of the $G_{l}$-s means that they have just one non-zero
entry in each row and column, which is either $+1$ or $-1$ by orthogonality. 

But this is not the end of the story, for one may show \cite{CMP}
that \begin{equation}
G_{l}^{-1}TG_{l}=\sigma_{l}^{2}\left(T\right)\,,\label{eq:Gcomm}\end{equation}
relating the Galois action on $S$ and $T$ via the matrices $G_{l}$.
This result has many important consequences, for example it implies
the upper bound, explained at the end of the previous section, for
the conductor in terms f the number of simple objects. Another important
consequence is that it allows to express \cite{CG2,Bauer} the matrices
$G_{l}$ in terms of $S$ and $T$:\begin{equation}
G_{l}=S^{-1}T^{l}ST^{m}ST^{l}\,,\label{eq:Gformula}\end{equation}
where $m$ denotes the inverse of $l$ in $\zn{N}$. In particular,
the group generated by the matrices $T$ and $S$ is invariant under
the Galois transformations $\sigma_{l}$.

\section{The projective kernel\label{sec:The-projective-kernel}}

Recall from section \ref{sec:The-modular-representation} that the
kernel of the modular representation is a congruence subgroup of level
$N$, i.e. it contains $\Gamma\left(N\right)$, where $N$ is the
order of the matrix $T$. It is natural to ask whether one could give
a more precise description of this subgroup. As it turns out, the
really interesting object is not the kernel itself, but the projective
kernel $\pk$, i.e. the subgroup of $\SL$ elements that are represented
by scalar multiples of the identity matrix, and for this one has a
nice description \cite{projker}. To explain it, we have to introduce
some more notions. First of all, let $K$ denote the projective order
of $T$, i.e. the smallest positive integer such that $T^{K}$ belongs
to the projective kernel $\pk$. Clearly, $K$ is a divisor of the
conductor $N$, and one may even show \cite{CMP} that the ratio $e=N/K$
is a divisor of $12$ (all divisors of $12$ arise this way, but for
a given value of $e$ one gets restrictions on $K$, e.g. $K$ has
to be odd if $4$ divides $e$). Next, consider \begin{equation}
\mathfrak{h}=\set{l\mod K}{G_{l}\in\mathbb{K}}\,,\label{eq:hdef}\end{equation}
which is a subgroup of $\zn{K}$ - actually it is a subgroup of exponent
$2$, as a corollary of Eq.(\ref{eq:Gcomm}). Note that we don't need
to know $\pk$ in order to determine $K$ and $\mathfrak{h}$, we
just have to check whether the relevant matrices are multiples of
the identity. To $\mathfrak{h}$ is associated the subgroup \begin{equation}
\Gamma\left(K,\mathfrak{h}\right)=\set{\left(\begin{array}{cc}
a & b\\
c & d\end{array}\right)\in\SL}{a,d\in\mathfrak{h}\,\,\mathrm{and}\,\,\Mod{b,c}{0}{K}}\:\label{eq:Gammadef}\end{equation}
 of $\SL$, and it is pretty easy to show that $\FA{\mathfrak{h}}$
is always a subgroup of the projective kernel $\pk$.

A more involved argument leads to the result that the index $\left[\pk:\FA{\mathfrak{h}}\right]$
is either $1$ or $2$. The common case is when the index is $1$,
i.e. $\pk=\FA{\mathfrak{h}}$. The index $2$ case is related to the
existence of so-called Galois currents \cite{projker}, which are
subtle symmetries of the modular representation, but is very rare
(in this case the ratio $e=N/K$ has to be odd, and $K$ a multiple
of $16$). But even in this last case one has canonical representatives
of the non-trivial $\FA{\mathfrak{h}}$ coset, i.e. the elements of
the projective kernel are known explicitly in terms of $K$ and the
subgroup $\mathfrak{h}$, which are easy to determine.

\section{Trace identities\label{sec:Trace-identities}}

Up to now we have been concerned with the properties of the modular
representation, i.e. the first term of the tower $\FB{}_{1},\FB{}_{2},\ldots$
of mapping class group modules. As we'll see in the present section,
this first term does determine the subsequent terms to a great extent.

The underlying idea is that a Modular Tensor Category associates to
each closed 3-manifold a number, called its partition function, which
is a topological invariant, i.e. it is the same for homeomorphic manifolds
\cite{Witten,Turaev}. While it is in general a difficult problem
to determine the value of this topological invariant for a given 3-manifold,
there are some classes of manifolds where the answer is known \cite{Freed,Rozansky}.

One such class of 3-manifolds are mapping tori. These arise through
the following construction: one starts with a closed surface $\FD{}$
of genus $g$, and one forms its Cartesian product with the unit interval
$\left[0,1\right]$. The resulting 3-manifold $\FD{}\times\left[0,1\right]$
is not closed, having two boundary components, each homeomorphic to
$\FD{}$. Identifying these two boundary components via some self-homeomorphism
$\phi$ of $\FD{}$ one obtains a closed 3-manifold called the mapping
torus of $\phi$. As it turns out, mapping tori of maps belonging
to the same mapping class are homeomorphic, i.e. the above construction
assigns to each mapping class in $\phi\in\FC{}_{g}$ a 3-manifold
$\FE{}{}$, which is well-defined up to homeomorphism. 

Given a Modular Tensor Category with associated sequence $\FB{}_{1},\ldots$
of mapping class group modules and a mapping class $\phi\in\FC{}_{g}$,
the partition function of the mapping torus $\FE{}{}$ is nothing
but the trace $\mathrm{Tr}\left(\phi\right)$ of $\phi$ on the corresponding
module $\FB{}_{g}$.

Another class of 3-manifolds of interest to us are the so-called Seifert-manifolds.
They have several equivalent characterizations, but the following
one seems the best suited for our purposes: we start again with a
closed surface $\FD{}$, but first we cut out from this surface $n>0$
non-overlapping disks, resulting in a surface $\FD{}^{*}$. We then
take the product of this surface with the circle $S^{1}$, which results
in a 3-manifold $\FD{}^{*}\times S^{1}$, which is not closed, its
boundary consisting of $n$ disjoint 2-tori. We obtain a closed 3-manifold
by pasting in solid tori to each of these boundary components, where
the pasting is characterized by a sequence $m_{1},\ldots,m_{n}$,
with each $m_{i}$ being a self-homeomorphism of the corresponding
2 dimensional torus. Once again, the homeomorphism type of the resulting
3-manifold does only depend on the mapping class of the $m_{i}$-s,
i.e. on the corresponding sequence of $\FC{}_{1}=\SL$ elements. In
this way we construct the Seifert-manifold $S\left(g;m_{1},\ldots,m_{n}\right)$,
where $g$ and $n$ are positive integers, while $m_{1},\ldots,m_{n}\in\SL$.
Using Dehn-surgery, one may show that the partition function of the
Seifert-manifold $S\left(g;m_{1},\ldots,m_{n}\right)$ equals \begin{equation}
\sum_{p}S_{0p}^{2-2g-n}\prod_{i=1}^{n}\left[m_{i}\right]_{0p}\,,\label{seifert}\end{equation}
where $\left[m_{i}\right]_{pq}$ stands for the matrix elements of
$m_{i}\in\FC{}_{1}=\SL$ in the distinguished basis of $\FB{}_{1}$. 

All this told, the point is to recognize that there is an overlap
between the class of Seifert-manifolds and mapping tori \cite{Jgeo}:
a mapping torus $\FE{}{}$ is homeomorphic to a suitable Seifert-manifold
$S\left(g^{*};m_{1},\ldots,m_{n}\right)$ if and only if the mapping
class $\phi$ has finite order \cite{Hempel}. Moreover, for a given
$\phi\in\FC{}_{g}$ of finite order, the integer $g^{*}$ and the
sequence $m_{1},\ldots,m_{n}\in\SL$ may be determined from the properties
of the branched cover $\FD{}\rightarrow\FD{}/\phi^{*}$, where $\phi^{*}$
denotes the lift of the mapping class $\phi$ to an automorphism of
the surface $\FD{}$. Eq.(\ref{seifert}) means that, for mapping
classes $\phi\in\FC{}_{g}$ of finite order, the trace of $\phi$
over the $\FC{}_{g}$ module $\FB{}_{g}$ may be expressed as an algebraic
combination of modular matrix elements, in other words the $\SL$
module $\FB{}_{1}$ does completely determine the restriction of the
representations $\FB{}_{g}$ to any finite subgroup of $\FC{}_{g}$.

In particular, for the identity mapping class of $\FC{}_{g}$, whose
trace equals obviously the dimension of $\FB{}_{g}$, one gets the
celebrated formula of Verlinde \cite{MS}\begin{equation}
\dim\FB{}_{g}=\sum_{p}S_{0p}^{2-2g}\,.\label{eq:ver}\end{equation}

Let's look now at the $\FC{}_{1}$ module $\FB{}_{1}$. Finite order
elements of $\FC{}_{1}=\SL$ are all conjugate to a power of either
$\left(\begin{array}{cc}
0 & -1\\
1 & 0\end{array}\right)$ or $\left(\begin{array}{cc}
0 & -1\\
1 & 1\end{array}\right)$. Performing the above sketched analysis, one arrives at the expression
\begin{equation}
\sum_{p}\frac{M\left(4\right)_{0p}M\left(4\right)_{0p}M\left(-2\right)_{0p}}{S_{0p}}\label{eq:Strace1}\end{equation}
for the trace of the operator representing $\left(\begin{array}{cc}
0 & -1\\
1 & 0\end{array}\right)$ on the space $\FB{}_{1}$, where $M\left(k\right)=S^{-1}T^{k}S$,
as before. But this operator is represented by the matrix $S$, whose
trace is just the sum of its diagonal elements, i.e.\begin{equation}
\mathrm{Tr}\, S=\sum_{p}S_{pp}\,.\label{eq:Strace2}\end{equation}
According to the above, the two expressions Eqs.(\ref{eq:Strace1})
and (\ref{eq:Strace2}) should be equal, providing a fairly nontrivial
consistency requirement on the matrix elements of the modular representation
\cite{Jgeo,BV}. A similar analysis applied to the other finite order
elements of $\SL$ yields more nontrivial restrictions.

\section{Summary}

The structure of a Modular Tensor Category, which formalizes the notion
of the Moore-Seiberg data of Rational Conformal Field Theory, is pretty
nice and intricate. We have tried to give an overview of some of the
results related to the tower of mapping class group modules. Most
of these results were concerned with the properties of the modular
representation, i.e. the first term of the tower. Already at this
level we can see some beautiful structures emerging, e.g. the Galois
action on modular matrices. Even more intriguing is the way the modular
representation influences the structure of the higher genus terms
of the tower, as explained in section \ref{sec:Trace-identities}. 

All these results seem to indicate that we are still far from the
end of the story, and that their are plenty of niceties to unravel.
A major issue is to find effective computational tools to deal with
Modular Tensor Categories and to enumerate them systematically. Needless
to say, there's still much work to be done to get a better understanding
of the subject.

\bigskip{}
\emph{Work supported by grants OTKA T037674, T043582 and T047041.}

\end{document}